\documentclass[10pt]{amsart}
\usepackage{a4,amsmath,amssymb,amsthm,eucal}
\usepackage{epsfig}
\newcommand\N{\mathbb{N}}
\newcommand\R{\mathbb{R}}

\newcommand{\deb}{\rightharpoonup}
\newcommand{\destar}{\overset{*}\deb}

\newcommand{\forget}[1]{}
\newcommand{\HH}{{\mathcal H}}

\def\Om{{\Omega}}  
\def\om2{{\Om\times\Om}}

\def\dist{\mathrm{dist}\,}

\newcommand{\lcal}{\mathcal{L}}

\newcommand{\pical}{\mathcal{P}}

\newcommand{\F}{\mathfrak{F}}

\newcommand\huno{{\mathcal H}^1}

\newcommand{\de}{\textrm{d}}
\newcommand{\res}{\llcorner} 

\newtheorem{theorem}{Theorem}[section]
\newtheorem{definition}[theorem]{Definition}

\newtheorem{lemma}[theorem]{Lemma}
\newtheorem{proposition}[theorem]{Proposition}
\newtheorem{corollary}[theorem]{Corollary}

\newtheorem{problem}[theorem]{Problem}

\theoremstyle{remark}
\newtheorem{remark}[theorem]{Remark}
\newtheorem{example}[theorem]{Example}

\numberwithin{equation}{section}

\author{A. Brancolini\and G. Buttazzo\and F. Santambrogio\and  E. Stepanov}
\thanks{The work of the fourth author has been partially supported by the INDAM project
``Traffic flows on complex networks'' and by GNAMPA}
\title[Long-term versus short-term planning]{Long-term planning versus short-term planning in the asymptotical location problem}

\begin{document}

\begin{abstract}
Given the probability measure $\nu$ over the given region
$\Omega\subset \R^n$, we consider the optimal location of a set
$\Sigma$ composed by $n$ points $\Om$ in order to minimize the
average distance $\Sigma\mapsto \int_\Om \dist(x,\Sigma)\,d\nu$ (the
classical optimal facility location problem). The paper compares two
strategies to find optimal configurations: the long-term one which
consists in
 placing all $n$ points at once in an optimal position, and the
 short-term one which consists in placing the points one by one adding
at each step at most one point and preserving the configuration
built at previous steps. We show that the respective optimization
problems exhibit qualitatively different asymptotic behavior as
 $n\to\infty$, although the optimization costs in both cases have the same asymptotic
 orders of vanishing.
\end{abstract}

\maketitle

\section{Introduction}\label{sec1}

Planning
an economic activity is in general an extremely complex
problem, where a high number of parameters often intervene. In
addition, the attitude of the planner has to be taken into account:
long-term planners take their decision through an optimization
process over a large time horizon, while short-term planners behave
optimizing day-by-day their strategies. Usually, the first kind of
behavior is perceived as more virtuous and efficient, while the
second is seen as easier to implement.

In this paper we analyze the long-term and short-term strategies in
a very simple model problem, and we give a way to measure the
efficiency of the first versus the second. The problem we consider
is the so-called {\em location problem\/} which can be roughly
described as follows. Suppose one has to open a certain given number
$n\in \N$ of identical facilities (e.g.\ plants, shops, distribution
centers, cinemas etc.) in the given urban region $\Omega\subset
\R^d$ which will be modeled by a compact convex set. If the density
of population in $\Omega$ is given by a known Borel probability
measure $\nu$, then the simplest way to measure how good is the
chosen configuration of facilities modeled by a set $\Sigma\subset
\Omega$ consisting of at most $n$ points (i.e.\ $\#\Sigma\leq n$),
is clearly to calculate the average distance the people have to
cover to reach the nearest facility
\begin{equation}\label{defiF}
F(\Sigma):=\int_\Omega\dist(x,\Sigma)d\nu(x).
\end{equation}
Hence the owner of the facilities is interested in locating them in
such a way as to minimize $F$. That is why, if he is able to open
all the facilities at once, he would choose a configuration
$\Sigma=\Sigma_n$ solving the following problem.

\begin{problem}\label{pb_long1}
Minimize the functional $\Sigma\mapsto F(\Sigma)$ subject to the
constraints
\[
\Sigma\subset\Omega,\,\#\Sigma\leq n.
\]
\end{problem}

We will further refer to such problem as a long-term planning
problem since it models the optimal choice of facility location so
as to satisfy global (long-term) needs of the facility owner. If
$\Sigma=\Sigma_n$ solves Problem~\ref{pb_long1}, we will denote
$l_n:=F(\Sigma_n)$ the respective optimal cost. It is worth
mentioning that such a problem has been extensively studied (see
e.g.~\cite{SuzDrez96,SuzOk95,MorgBolt01} for recent surveys on the
subject), but nevertheless a lot of interesting and important
questions regarding this problem still remain open.

If however the owner of the facilities is unable to open all the
facilities immediately (say, if he does not have enough financial
resources to do that), he will try to open them one by one, trying
to minimize the average distance functional $F$ at each step (i.e.\
when opening each facility), but, of course, taking into
consideration the location of facilities already opened at previous
steps. Such a short-term optimization strategy amounts to solving
the following problem (in the sequel referred to as short-term
planning problem).

\begin{problem}\label{pb_short1}
For each $n\in \N$, $n\geq 1$, find a set $\Sigma=\Sigma_n'$
minimizing the functional $\Sigma\mapsto F(\Sigma)$ subject to the
constraints
\[
\Sigma'_{n-1}\subset\Sigma\subset\Omega,\, \#\Sigma\leq n,
\]
where $\Sigma'_0:=\emptyset$.
\end{problem}

The above model, to the best of our knowledge, can be traced back
to~\cite{TSuzAsamOkab91}, where the short-term facility allocation
strategy is called a {\em myopic allocation policy\/} (see also
references therein for similar formulations). In the same paper the
authors propose also several other allocation strategies which can
be considered as intermediate between the short-term and the
long-term ones.

We will further denote $s_n:=F(\Sigma_n')$ the optimal cost of the
above short-term optimization problem.

In this paper we discuss the asymptotic
behavior of solutions to the above two problems as $n\to \infty$, by
studying the weak limits (in a suitable sense) of optimal
configurations $\Sigma_n$ and $\Sigma'_n$ as well as the optimal
costs $l_n$ and $s_n$. Namely, we will be interested in finding
answers to the following  questions.
\begin{itemize}
\item[(A)] Find the asymptotic order of $l_n$ (resp.\ $s_n$)
as $n\to \infty$, i.e.\ find an exponent $\alpha>0$ such that
\[
C_1 n^{-\alpha} \leq l_n \leq C_2 n^{-\alpha} \qquad\qquad
(\mbox{resp. } C_1 n^{-\alpha} \leq s_n \leq C_2 n^{-\alpha})
\]
for some positive constants $C_1$ and $C_2$ and for $n$ sufficiently
large. To simplify the notation, in the sequel we will write in this
case $l_n\sim n^{-\alpha}$ (resp.\ $s_n\sim n^{-\alpha}$).
\item[(B)] Find precise asymptotic estimates for $l_n$ (resp.\ $s_n$), i.e.\ find $\lim_n n^{\alpha}l_n$
(resp.\ $\lim_n n^{\alpha} s_n$), or just $\liminf$ and $\limsup$,
should the limit not exist.
\item[(C)] Describe the asymptotic behavior of the minimizers, i.e. find all the weak limits,
in a suitable sense, of subsequences of minimizers.
\end{itemize}

In the pioneering paper~\cite{TSuzAsamOkab91} some accurate
numerical calculations of optimal short-term configurations for the
case of the uniform density on a line and on a two-dimensional
square have been provided for rather small values of $n$. Moreover,
some deep insights have been formulated regarding the behavior of
solutions. Nevertheless, the above questions have not been
explicitly addressed nor even formally posed. In this paper we
provide rigorous results concerning the nature of the problem which
partially confirm the insights of~\cite{TSuzAsamOkab91}, and also go
further. In Section~\ref{sec2} we summarize all the known results on
the asymptotic behavior of solutions to the long-term optimization
Problem~\ref{pb_long1} that we need for the purpose of comparison
with the short-term optimization Problem~\ref{pb_short1}. In
Sections~\ref{sec3} and~\ref{onedim} we show that although the
optimal cost of the short-term optimization Problem~\ref{pb_short1}
has the same order of asymptotic expansion as $n\to \infty$ (i.e.\
the answer to question~(A) is the same for both problems),
nevertheless, the two problems are qualitatively different even in
the simplest one-dimensional case $d=1$ (i.e.\ the answers to
questions~(B) and~(C) are qualitatively different for these
problems). Finally, we conclude the paper by some remarks and open
questions supported by numerical evidence, as well as the
description of a similar problem we feel important for applications.

\section{The long-term problem}\label{sec2}

The asymptotic behavior of the long-term optimal location
Problem~\ref{pb_long1} has been intensively studied both using
geometric (see e.g.~\cite{FejesToth} and~\cite{MorgBolt01}) and
variational methods~\cite{BouJimRaj,TillMosc02}, in the latter case
mainly by means of $\Gamma-$convergence tools
(see~\cite{introgammaconve} for details on the theory). We summarize
here the most important properties of this problem, in order to make
later a comparison with the respective properties of the short-term
problem.

Notice first that
\begin{equation}\label{long,n-1/d}
l_n \leq Cn^{-1/d}.
\end{equation}
This estimate is straightforward, if one considers a set $\Sigma$
composed by $n$ points placed on a uniform grid of size
approximately equal to $n^{-1/d}$.

We use $\Gamma-$convergence theory to find the answers to the
questions~(A)--(C) posed in the Introduction. In order to apply it,
we need to work with functionals defined on a common space, which we
choose to be the space of all Borel probability measures
$\pical(\Omega)$. To this aim we identify each set
$\Sigma\subset\Omega$ having $\#\Sigma<+\infty$ with the measure
$\mu_{\Sigma}\in\pical(\Omega)$ defined by
\[
\mu_{\Sigma}:=\frac{1}{\#\Sigma}\sum_{x\in\Sigma}\delta_{x}.
\]
We
define now a sequence of functionals on the space $\pical(\Omega)$
by setting
\begin{equation}\label{eq_FnProb}
\F_n(\mu):=\begin{cases}n^{1/d}F(\Sigma),&\text{ if }\mu=\mu_{\Sigma},\,\#\Sigma\leq n,\\
                        +\infty ,          &\text{otherwise.}\end{cases}
\end{equation}
The coefficient $n^{1/d}$ in the above formula prevents the
minimization from degenerating and is chosen according
to~\eqref{long,n-1/d}. It is straightforward to recognize that
minimizing $\F_n$ is equivalent to solving Problem~\ref{pb_long1} up
to the above identification of sets with probability measures.

We give here the $\Gamma-$convergence result only for the case
$\nu=f\cdot{\mathcal L}^d$. This result could easily be generalized
to a generic measure $\nu$ by inserting in the $\Gamma-$limit only
the absolutely continuous part of $\nu$ with respect to $\lcal^d$.
The proof of this theorem, for the case when $f$ is a lower
semicontinuous function, can be found in~\cite{BouJimRaj}. For the
proof in the general case we provide below (i.e.\ when $f\in
L^1(\Omega)$ is not necessarily lower semicontinuous), one can apply
the technique developed in~\cite{TillMosc02} for the study of a
similar problem (the so-called {\it irrigation problem}).

\begin{theorem}\label{th_UCGamloc}
The sequence of functionals $\{\F_n\}_n$ $\Gamma$-converges with
respect to the weak$^*$ convergence of measures to the functional
$\F_\infty$: ${\mathcal P}(\Omega)\to \bar{\R}$ defined by the
formula
$$
\F_\infty(\mu):=
\begin{array}{rl}
\displaystyle
\theta_{d} \int_\Omega \frac{f}{\rho^{1/d}}\,d\lcal^d,
& \mbox{ where }\mu=\rho\cdot{\mathcal L}^d + \mu_{sing},
\end{array}
$$
where $\rho\in L^1(\Omega)$, $\mu_{sing}$ stands for the singular
part of the measure $\mu$ with respect to $\lcal^d$, and
 $\theta_d$ is a constant depending only on the
dimension $d$ satisfying $0<\theta_d<\infty$ and given by
$$\theta_d:=\inf\left\{\liminf_n n^{1/d}\int_{[0,1]^d}
\dist(x,\Sigma_n)dx\;:\;\Sigma_n\subset [0,1]^d,\,\#\Sigma_n\leq
n\right\}.$$
\end{theorem}

The most important consequences of this theorem are summarized in
the following corollary and answer the questions we are interested
in.

\begin{corollary}\label{consequences}
The following assertions hold.
\begin{itemize}
\item[(A),(B)]
One has
\[ \lim_n n^{1/d} l_n =\min
\big\{\F_{\infty}(\mu)\,:\,\mu\in\pical(\Omega)\big\}=\theta_d||f||_{d/(d+1)}>0,
\]
hence,  in particular, $l_n\sim n^{-1/d}$.
\item[(C)] Denoting by $\mu_n:=\mu_{\Sigma_n}$ the measures associated to a sequence of minimizers for
Problem~\ref{pb_long1}, we have  that $\mu_n\destar\bar{\mu}$ as
$\nu\to\infty$ in the weak$^*$ sense of measures, where $\bar{\mu}$
is the unique minimizer of $\F_{\infty}$ and is given by the formula
\[
\bar{\mu}=cf^{d/(d+1)}\cdot\lcal^d\mbox{ with } c:=\left(\int_\Omega
f^{d/(d+1)}\,d\lcal^d\right)^{-1}.
\]
 In particular, if $\nu$ has constant
density, i.e.\ $\nu=c\cdot\lcal^d$, then $\bar{\mu}$ has constant
density as well, namely, $\bar{\mu}=\nu$.
\end{itemize}
\end{corollary}

\begin{proof}
This statement comes from well-known properties of
$\Gamma-$convergence (i.e. convergence of minima and of minimizers),
once we find the unique minimizer of $\F_{\infty}$. Due to the
strictly decreasing nature of $\F_{\infty}$ with respect to  $\rho$
on $\{f>0\}$ and to the fact that $\mu_{sing}$ and $\rho
1_{\{f=0\}}$ do not affect the value of $\F_{\infty}$, it is
straightforward that the minimizers should be absolutely continuous
and concentrated on $\{f>0\}$. To identify the density of the
absolutely continuous part, set
$\lambda:=f^{d/(d+1)}\cdot\lcal^d\res \Omega$, and $w:=\rho
f^{-d/(d+1)}$. With this notation, finding the minimizers to
$\F_{\infty}$ is equivalent to minimizing the functional $w\mapsto
\int_{\Omega}w^{-1/d}\,d\lambda $ over the set
$$\left\{w\in L^1(\Omega,\lambda),\,w\geq 0,\,\int_{\Omega}w\,d\lambda=1\right\} .$$
By convexity of the map $w\mapsto w^{-1/d}$, it immediately follows
from Jensen inequality that the minimum of the latter functional is
attained at a constant function $w$. This shows
$\bar{\mu}=cf^{d/(d+1)}\cdot\lcal^d$ and the computation of the
constant $c$ follows from the constraint
$\bar{\mu}\in\pical(\Omega)$. At last, the value of
$\min\F_{\infty}$ is obtained by plugging $\bar{\mu}$ into the
expression for $\F_{\infty}$.
\end{proof}

The exact values of the constants $\theta_d$ are known in the
one-dimensional and two-dimensional cases. Namely, if $\Omega=[0,1]$
and $\nu=\lcal^{1}\res \Omega$ it is actually easy to
 compute explicitly the unique minimizer of Problem~\ref{pb_long1} which
 is given by the set
of $n$ points located at the centers of $n$ equal disjoint intervals
forming a partition of $\Omega$. In other words, one has
$$\Sigma_n=\bigcup_{i=1}^n\left\{\frac{2i-1}{2n}\right\},\quad F(\Sigma_n)=\frac{1}{4n},\quad\mbox{ so that }\quad\theta_1=\frac{1}{4}.$$
In the two-dimensional case when $\Omega=[0,1]^2$ and
$\nu=\lcal^{2}\res \Omega$ it is known that the configuration of $n$
points placed in centers of regular hexagons, is asymptotically
optimal as $n\to \infty$ (see~\cite{FejesToth}
or~\cite{MorgBolt01}), which gives possibility to compute explicitly
the constant $\theta_2$. Namely one gets
$$\theta_{2}=\int_\sigma |x|\, dx=\frac{3\log 3 +4}{6\sqrt{2}\,3^{3/4}}\approx 0.377$$
where $\sigma\subset \R^2$ stands for the regular hexagon of unit
area centered at the origin.

\section{The short-term problem}\label{sec3}

In this section we answer question~(A) posed in the Introduction
regarding the short-term optimal location Problem~\ref{pb_short1}.
Namely, we will show that
\begin{equation}\label{eq_snapprox1}
s_n\sim n^{-1/d}
\end{equation}
whenever $\nu\ll \lcal^d$,
 similarly to the asymptotic estimate $l_n\sim n^{-1/d}$ proved
 in  Section~\ref{sec2}.

Before proving~\eqref{eq_snapprox1} we find it important to remark
that the values $s_n$ actually may depend not only on $n$, but also
on the chosen sequence of solutions $\{\Sigma'_n\}$ of the
short-term Problem~\ref{pb_short1}. This is due to the fact that at
each minimization step both the position of the next point and the
new minimum value may depend on the history, i.e.\ on the
configuration chosen on the previous steps. In other words, the
choice of the optimal set at each step may affect the minimal values
of all the following steps  as the following example shows.

\begin{example}\label{ex_nonuniquesn}
Consider the one-dimensional situation $d=1$ with $\Om:=[0,4]$ and
\[
\nu:=\left(\frac 1 2 \cdot 1_{[0,1]} + \frac 1 4
1_{[2,4]}\right)\cdot \lcal^1.
\]
 Then one clearly
has that both singletons $\{1\}$ and $\{2\}$ are solutions to the
short-term location problem at the first step $n=1$, and both give
the value
 $s_1=5/4$. Now, if one takes $\Sigma'_1:=\{1\}$, then
at the second optimization step $n=2$ we get the unique minimizer
$\Sigma'_2:=\{1,3\}$, which gives the value $s_2=1/2$. On the other
hand, if at the first step one takes $\Sigma'_1:=\{2\}$, then at the
second step one gets the unique minimizer $\Sigma'_2:=\{1/2,2\}$
which gives a different value $s_2=5/8$.
\end{example}

Coming back to proving~\eqref{eq_snapprox1}, we observe that it is
impossible to use the $\Gamma$-convergence theory for this purpose.
Namely, the constraint $\Sigma'_{n+1}\supset\Sigma'_n$ which is
imposed at each minimization step, actually gives raise to a
sequence of problems that, once rescaled and expressed in terms of
probability measures as in Section~\ref{sec2}, involve the
functionals $\F'_n$ given by
$$\F'_n(\mu)=
\left\{
\begin{array}{rrl}\F_n(\mu),&\mbox{if
  either }  &\mu=\frac{n-1}{n}\mu_{\Sigma'_{n-1}}+\frac{1}{n}\delta_{x},\,x\in\Omega, \\
    & \mbox{or } &\mu=\mu_{\Sigma'_{n-1}},
\\
                         +\infty, &\mbox{ otherwise.}
\end{array}
\right.
$$
It is therefore not difficult to see that whenever
$\mu_{\Sigma'_n}\destar\bar{\mu}$ in weak$^*$ sense of measures as
$n\to\infty$, then the $\Gamma-$limit functional $\F'_{\infty}$
would be finite only on $\bar{\mu}$ itself. This means that taking
the limit of a sequence $\mu_{\Sigma_n}$ and using the fact that it
minimizes $\F'_{\infty}$ gives no additional information on the
limit itself. Hence it is not possible to find $\bar{\mu}$ in this
way, contrary to the long-term case. We will therefore analyze the
short-term Problem~\ref{pb_short1} directly, without using
$\Gamma-$convergence tools.

The following assertion is valid.

\begin{theorem}\label{shortth}
Let 
$\nu\ll \lcal^d$. Then
\[
C_1 n^{-1/d} \leq s_n \leq C_2 n^{-1/d},
\]
where the positive constants $C_1$ and $C_2$ do not depend on the
choice of the sequence of solutions to the short-term optimal
location Problem~\ref{pb_short1}. In
particular,~\eqref{eq_snapprox1} holds. 
\end{theorem}

\begin{proof}
Let $\{\Sigma'_n\}$ (resp.\ $\{\Sigma_n\}$) be a sequence of
minimizers for the short-term optimal location
Problem~\ref{pb_short1} (resp.\ long-term optimal location
Problem~\ref{pb_long1}), where $\nu=f\cdot\lcal^d$, $f\in
L^1(\Omega)$. Of course, the long-term cost is lower than the
short-term one, namely,
$$F(\Sigma_n)\le F(\Sigma'_n).$$
This provides the required estimate from below as a consequence of
Corollary~\ref{consequences}. It remains to prove  an estimate of
the form
$$F(\Sigma'_n)\le Bn^{-1/d}$$
for a suitable constant $B$ independent of the choice of the
sequence of solutions to the short-term optimal location
Problem~\ref{pb_short1}.

Once we have a set $\Sigma'_n$ we want to estimate how much the
functional $F$ decreases when we add a point $x_0\in\Omega$. For the
sake of brevity denote $\delta(x):=\dist(x,\Sigma'_n)$. If we set
$\Sigma:=\Sigma'_n\cup\{x_0\}$, it is clear that for $x\in
B(x_0,\delta(x_0)/4)$ we have $\dist(x,\Sigma)<\delta(x_0)/4$,
$\delta(x)>\frac{3}{4}\delta(x_0)$ and hence
$\dist(x,\Sigma)<\delta(x)-\delta(x_0)/2$. Thus, if we set
\[
g(x_0):=\nu\left(B\left(x_0,\frac{\delta(x_0)}{4}\right)\right)\frac{\delta(x_0)}{2},
\]
 we get
$$F(\Sigma)\leq F(\Sigma'_n)-g(x_0)\quad\mbox{  and  }\quad F(\Sigma'_{n+1})\leq F(\Sigma'_n)-\sup_{x_0\in\Omega} g(x_0).$$
To estimate $\sup_{x_0\in\Omega} g(x_0)$ we use the inequalities
\begin{align*}
\int_{\Omega}g(x_0)dx_0& =\int_{\Omega}dx_0\int_{\Omega}d\nu(x) \frac{\delta(x_0)}{2}1_{\{|x-x_0|<\delta(x_0)/4\}}(x)\\
&
=\int_{\Omega}d\nu(x)\int_{\Omega}\frac{\delta(x_0)}{2}1_{\{|x-x_0|<\delta(x_0)/4\}}(x_0)dx_0\\
& \geq
\int_{\Omega}d\nu(x)\int_{\Omega}\frac{\delta(x_0)}{2}1_{\{|x-x_0|<\delta(x)/5\}}(x_0)dx_0\\
&\geq\int_{\Omega}d\nu(x)\frac{4\omega_d}{2\cdot
5^{d+1}}\delta(x)^{d+1}
\geq\frac{2\omega_d}{5^{d+1}}\left(\int_{\Omega}d\nu(x)\delta(x)\right)^{d+1},
\end{align*}
where $\omega_d$ stands for the volume of the unit ball in $\R^d$.
In the above chain of inequalities we have used that, for $x_0\in
B(x,\delta(x)/5)$ the condition $|x-x_0|<\delta(x_0)/4$ is always
satisfied as well as the inequality $\delta(x_0)\geq
\frac{4}{5}\delta(x)$; moreover, the last inequality in the above
chain is an application of Jensen inequality since $\nu$ is a
probability measure.

We have therefore obtained
$$\sup_{x_0\in\Omega} g(x_0)\geq \frac{1}{|\Omega|}\int_{\Omega}g(x_0)dx_0\geq C\left(\int_{\Omega}d\nu(x)\delta(x)\right)^{d+1}
=CF(\Sigma'_n)^{d+1},$$ which implies
$$F(\Sigma'_{n+1})\leq F(\Sigma'_n)-CF(\Sigma'_n)^{d+1},$$
where $C>0$ depends only on the dimension $d$ of the underlying
space.
 The conclusion follows now from Lemma~\ref{recursive_inequality} below
(minding that the value $s_1$ is independent of the choice of the
sequence of solutions to the short-term optimal location
Problem~\ref{pb_short1}).
\end{proof}

\begin{lemma}\label{recursive_inequality}
Let $\{a_n\}$ be a sequence of nonnegative numbers satisfying
$a_{n+1}\leq a_n-Ca_n^{d+1}$ for all $n\in \N$, where $C > 0$. Then
there exists a number $B>0$ (depending only on $a_1$, $C$ and $d$),
such that $a_n\leq Bn^{-1/d}$ for all $n\in \N$.
\end{lemma}

\begin{proof}
The proof will be performed by induction, simultaneously with the
choice of $B$.

The step $n = 1$ is satisfied choosing $B \geq a_1$. Now we look for
a condition on $B$ such that the following statement is satisfied:
if for some $n\in \N$ one has $a_n \leq Bn^{-1/d}$, then $a_{n + 1}
\leq B(n + 1)^{-1/d}$. Since
$$
a_{n + 1} \leq a_n - Ca_n^{d + 1} \leq Bn^{-1/d} - Ca_n^{d + 1},
$$
then  we must impose that
$$
Bn^{-1/d} - Ca_n^{d + 1} \leq B(n + 1)^{-1/d}.
$$
Using again the fact that $a_n\leq Bn^{-1/d}$, we get
$$
Bn^{-1/d} \leq B(n + 1)^{-1/d} + CB^{d + 1}n^{-(d + 1)/d},
$$
or equivalently
$$
\frac{\left(1 + \frac{1}{n}\right)^{1/d} - 1}{1/n} \leq CB^d\left(1
+ \frac{1}{n}\right)^{1/d}.
$$
Minding that
$$
\frac{\left(1 + \frac{1}{n}\right)^{1/d} - 1}{1/n} \nearrow
\frac{\de(x^{1/d})}{\de x}|_{x = 1} = 1/d,
$$
it is sufficient to impose $1/d \leq CB^d\left(1 +
1/n\right)^{1/d}$, or even $1/d\leq CB^d$. Therefore, the choice
$B:= a_1\vee (1/dC)^{1/d}$ suffices for the proof to be concluded by
induction.
\end{proof}


\section{The one-dimensional case with uniform measure}\label{onedim}

In this section we address questions~(B) and~(C) posed in the
Introduction regarding the short-term optimal location
Problem~\ref{pb_short1}. In particular, studying question~(B), we
will show that, even in a very simple one-dimensional situation with
$\nu$ being the uniform measure on an interval, the ratio $s_n/l_n$
(which is never smaller than one) does not tend to a limit as $n\to
\infty$, and hence, neither does $n^{1/d} s_n$ (although by
Theorem~\ref{shortth} one has $n^{1/d} s_n\sim 1$), since by
Corollary~\ref{consequences} $\lim_n n^{1/d} l_n$ exists. Moreover,
we will show that in this case even $\liminf_n s_n/l_n
>1$.
We further show that in the same situation there are infinitely many
limit measures (in the weak$^*$ sense) of the sequences of
$\mu_{\Sigma'_n}$ as $n\to \infty$, where $\Sigma'_n$ solves
Problem~\ref{pb_short1}, and what is more, neither of such limit
measures is equal to the unique limit measure of the sequence of
solutions to the long-term Problem~\ref{pb_long1}.

In this section we restrict ourselves to the case $\Om=[0,1]$ and
$\nu=\lcal^1$. We will further identify each set $\Sigma\subset \Om$
having finite number of points with the partition of $\Om$ by the
points of $\Sigma$. Namely, if $\Sigma=\{x_1,\ldots, x_k\}$, we will
always order the elements of $\Sigma$ in such a way that
\[
x_1\leq x_2\leq\ldots\leq x_k,
\]
and identify $\Sigma$ with the partition of $\Om$ given by the
intervals $(\Delta_i)_{i=1}^{k+1}$, where
\[
\Delta_1:=[0, x_1], \Delta_2:=[x_1,x_2],\ldots, \Delta_k=[x_{k-1},
x_k], \Delta_{k+1}:=[x_k, 1].
\]
The intervals $\Delta_1$ and $\Delta_{k+1}$ will be further called
\emph{external intervals}, while all the other intervals of this
partition will be called \emph{internal intervals}. We will also
identify the same set
 with the
$(k+1)$-dimensional vector whose entries are the lengths of the
intervals of the respective partition
$$(|\Delta_1|,\ldots,|\Delta_{k+1}|)=(x_1, x_2-x_1, x_3-x_2,\ldots, x_k-x_{k-1}, 1-x_k),$$
and we write
\[
\Sigma\simeq (x_1, x_2-x_1, x_3-x_2,\ldots, x_k-x_{k-1}, 1-x_k).
\]
 For
instance, the set $\{1/3, 2/3\}$ is identified with the respective
partition of $[0,1]$, and with the vector $(1/3, 1/3, 1/3)$, i.e.\
$\{1/3, 2/3\}\simeq (1/3, 1/3, 1/3)$. The intervals $[0, 1/3]$ and
$[2/3, 1]$ are external intervals while the interval $[1/3, 2/3]$ is
internal.

\begin{proposition}\label{prop_divis1}
Given a $\Sigma_k\subset \Om$, assume that $\Sigma_k\simeq (\lambda_1, \ldots, \lambda_{k+1}),$ and
let $\Sigma_{k+1}$ be a minimizer of $F$ over all sets
$\Sigma\subset \Om$ such that $\Sigma\supset \Sigma_k$ and
 $\#\Sigma_{k+1}=k+1$. Then  $\Sigma_{k+1}=\Sigma_k\cup
\{x\}$, while either of the two conditions hold.
\begin{enumerate}
  \item Either $x$ is the center of some internal interval
  $\Delta_i$, so that
  \[
\Sigma_{k+1}\simeq (\lambda_1, \ldots,\lambda_{i-1}, \lambda_i/2,
\lambda_i/2, \lambda_{i+1},\ldots, \lambda_{k+1}),
  \]
  while
  \[
F(\Sigma_{k+1})= F(\Sigma_k)- \lambda_i^2/8,
  \]
  \item or $x$ divides some external interval
  $\Delta_1$ or $\Delta_{k+1}$ with the length ratio $1:2$ closer to the
  boundary of $[0,1]$, so that
\[
\Sigma_{k+1}\simeq (\lambda_1/3, 2\lambda_1/3, \lambda_2, \ldots,
\lambda_{k+1}) \mbox{ or } \Sigma_{k+1}\simeq (\lambda_1\ldots,
\lambda_k, 2\lambda_{k+1}/3, \lambda_{k+1}/3)
  \]
  while
  \[
F(\Sigma_{k+1})= F(\Sigma_k)- \lambda_i^2/3,
  \]
  where $i=1$ (if $x\in \Delta_1$) or $i=k+1$ (if $x\in \Delta_{k+1}$).
\end{enumerate}
\end{proposition}

\begin{proof}
If $x\in \Delta_i$, $i\not\in \{1, k+1\}$ (i.e.\ $\Delta_i$ is
internal), then
\begin{align*}
F(\Sigma_{k+1}) &=
\int_{[0,1]\setminus\Delta_i}\dist(z,\Sigma_k)\, dz +
\int_{\Delta_i}\dist(z,\{x_{i-1}, x, x_i\})\, dz\\
& =\int_0^1 \dist(z,\Sigma_k)\, dz - \frac{(x_i-x)(x-x_{i-1})}{2},
\end{align*}
and to conclude, it is enough to note that the function
$$[x_{i-1}, x_i]\ni x\mapsto (x_i-x)(x-x_{i-1})/2$$
attains its maximum value $\lambda_i^2/8$ at $x=(x_{i-1}+x_i)/2$. The case when
$\Delta_i$ is external, say, $i=1$ (the case $i=k+1$ is completely
symmetric), is absolutely analogous, once one notes that
\begin{align*}
F(\Sigma_{k+1}) &=
\int_{[0,1]\setminus\Delta_1}\dist(z,\Sigma_k)\, dz +
\int_{\Delta_i}\dist(z,\{ x, x_1\})\, dz\\
& =\int_0^1 \dist(z,\Sigma_k)\, dz -
\frac{x_1^2}{2}-\left(\frac{x^2}{2}+\frac{(x_1-x)^2}{4}\right)\\
&= F(\Sigma_k)-\frac{3x^2-2x_1x+x_1^2}{4}.
\end{align*}
Indeed, one has that the function $[0,x_1]\ni x \mapsto
(3x^2-2x_1x+x_1^2)/4$ attains its maximum value $\lambda_1^2/3$
at $x=x_1/3$.
\end{proof}

We will further say that the point $x\in \Delta_i$ is in optimal
position, if either $\Delta_i$ is an internal interval and $x$ is
its center, or $\Delta_i$ is an external interval, and $x$ divides
it with the length ratio $1:2$ closer to the
  boundary of $[0,1]$. The above proposition says that
  whenever $\Sigma_k\subset \Om$, $\#\Sigma_k=k$ and $\Sigma_{k+1}\subset \Om$, $\#\Sigma_{k+1}=k+1$
  solve
Problem~\ref{pb_short1}, then $\Sigma_{k+1}=\Sigma_k\cup \{x\}$ with
$x$ in optimal position.

We set now
\[
\Om_i:= \left[\frac 1 {2\cdot 3^{i+1}},\frac 1 {2\cdot 3^{i}}
\right]\cup \left[1- \frac 1 {2\cdot 3^{i}},1-\frac 1 {2\cdot
3^{i+1}} \right], \qquad\qquad i\in \N,
\]
so that, clearly $\{\Om_i\}_i$ gives a partition of $\Om$, while
$|\Om_i| = 2/3^{i+1}$. This allows us to formulate the following
corollary to the above Proposition~\ref{prop_divis1}.

\begin{corollary}\label{col_divis1a}
Let $\Sigma'_k\subset \Om$, $\#\Sigma'_k=k$ be a solution to
Problem~\ref{pb_short1}. Then, for the corresponding partition of
$\Om$ one has that for each interval $\Delta_i$ there exists a
unique couple of
  numbers $(j,h)\in \N^2$, such that $|\Delta_i|=3^{-j}2^{-h}$,
  while
\begin{enumerate}
\item
  if $\Delta_i$ is internal, then $j\geq 1$ and
  $\Delta_i\subset \Om_j$.
  \item if $\Delta_i$ is external (i.e.\ $i=1$ or $i=k+1$), then
  $h=1$ and
\[
j=\left\{
\begin{array}{rl}
  0, & k=1, \\
  j_0\mbox{ or } j_0-1,
 & \mbox{otherwise},
\end{array}\right.
\]
where $j_0:=\sup\{j\,:\, \Om_{j}\supset \Delta_m\mbox{ for some }
m=2,\ldots, k\}$.
\end{enumerate}
\end{corollary}

\begin{proof} The proof is easily obtained by induction on $k$.
\end{proof}

Consider  an arbitrary $\Sigma'_k\simeq (\lambda_1, \ldots,
\lambda_{k+1})$ solving Problem~\ref{pb_short1}.
 Thanks to the corollary~\ref{col_divis1a}
we will identify each number $\lambda_i$ corresponding to the
internal interval $\Delta_i$ with the respective couple $(j,h)\in
\N^2$, $j\geq 1$, such that $\lambda_i=3^{-j}2^{-h}$. We will
further identify every external interval $\Delta_i$ with the couple
$(j+1/2,-1/2)$, where $j\in \N$ is such that $\lambda_i=3^{-j}
2^{-1}$ (the reason for the latter identification will be explained
in a moment).

Consider now the set $D$ of couples $(j,h)\in \R^2$, where $(j,h)\in
\N^2$, $j\geq 1$, or $j\in \N+1/2$, $h=-1/2$. We have that each
interval $\Delta_i$, and each $\lambda_i$ is identified with a
unique point $d_{\lambda_i}\in D$. We introduce the ordering on $D$
according to the following definition.

\begin{definition}\label{def_orderD}
We denote $(j_2, h_2)\succ (j_1, h_1)$, if
\[
h_2-h_1 + (j_2-j_1)\log 3/\log 2 <0.
\]
\end{definition}

Observe that the above relation well orders the set $D$. We have now
the following simple assertion.

\begin{proposition}\label{prop_divis2}
Let $\Sigma'_k\subset \Om$, $\#\Sigma'_k=k$ be a solution to
Problem~\ref{pb_short1}. Let $p, q\in \{1, \ldots, k+1\}$, and
consider $d_{\lambda_p}=(j_{\lambda_p}, h_{\lambda_p}),
d_{\lambda_q}=(j_{\lambda_q}, h_{\lambda_q})\in D$. Consider
$\Sigma_{k+1}^p:=\Sigma'_k\cup \{x^p\}$ and
$\Sigma_{k+1}^q:=\Sigma'_k\cup \{x^q\}$, where $x^p\in\Delta_p$ and
$x^q\in\Delta_q$ are in optimal positions. If $d_{\lambda_q}\succ
d_{\lambda_p}$, then
\[
F(\Sigma_{k+1}^p)> F(\Sigma_{k+1}^q).
\]
\end{proposition}

\begin{proof}
One considers separately three cases.

{\sc Case 1.} Both $\Delta_p$ and $\Delta_q$ are internal intervals.
Then
\begin{align*}
F(\Sigma_{k+1}^p) &= F(\Sigma'_k)- \lambda_p^2/8 \\
                          &= F(\Sigma'_k)- 3^{-2j_{\lambda_p}}
                          2^{-2h_{\lambda_p}}/8 \\
                          & > F(\Sigma'_k)- 3^{-2j_{\lambda_q}}
                          2^{-2h_{\lambda_q}}/8 = F(\Sigma_{k+1}^p),
\end{align*}
because $h_{\lambda_q}-h_{\lambda_p} +
(j_{\lambda_q}-j_{\lambda_p})\log 3/\log 2 <0$ implies
$3^{-2j_{\lambda_p}}
                          2^{-2h_{\lambda_p}} <3^{-2j_{\lambda_q}}
                          2^{-2h_{\lambda_q}}$.

{\sc Case 2.} One of the intervals (say, $\Delta_p$) is internal,
another one ($\Delta_q$) is external. Then, minding that
$h_{\lambda_p}=-1/2$, and hence $\lambda_p=3^{-j_{\lambda_p}+1/2}
2^{-1}=3^{-j_{\lambda_p}-1/2}2^{2h_{\lambda_p}}$, we get
\begin{align*}
F(\Sigma_{k+1}^p) &= F(\Sigma'_k)- \lambda_p^2/3 \\
                          &= F(\Sigma'_k)- 3^{-2j_{\lambda_p}}
                          2^{4h_{\lambda_p}} \\
                          & =F(\Sigma'_k)- 3^{-2j_{\lambda_p}}
                          2^{4h_{\lambda_p}}\\
                          & > F(\Sigma'_k)- 3^{-2j_{\lambda_q}}
                          2^{-2h_{\lambda_q}}/8 = F(\Sigma_{k+1}^p),
\end{align*}
because $h_{\lambda_q}-h_{\lambda_p} +
(j_{\lambda_q}-j_{\lambda_p})\log 3/\log 2 <0$ and
$h_{\lambda_p}=-1/2$ implies $3^{-2j_{\lambda_p}}
                          2^{-h_{\lambda_p}} <3^{-2j_{\lambda_q}}
                          2^{-2h_{\lambda_q}}/8$.

{\sc Case 3.} Both intervals $\Delta_p$ and $\Delta_q$ are external,
so that $\lambda_p=3^{-j_{\lambda_p}+1/2} 2^{-1}$ and
$\lambda_q=3^{-j_{\lambda_q}+1/2} 2^{-1}$, while $j_{\lambda_q} <
j_{\lambda_p}$.
 Then
\begin{align*}
F(\Sigma_{k+1}^p) &= F(\Sigma'_k)- \lambda_{\lambda_p}^2/3 \\
                          &= F(\Sigma'_k)- 3^{-2j_{\lambda_p}+1}
                          2^{-2}/3 \\
                          & > F(\Sigma'_k)- 3^{-2j_{\lambda_q}+1}
                          2^{-2}/3  = F(\Sigma_{k+1}^p),\qedhere
\end{align*}
\end{proof}

 As an immediate consequence of
Proposition~\ref{prop_divis2} we get the following theorem.

\begin{theorem}\label{th_div1}
Let $\Sigma'_k\subset \Om$, $\#\Sigma'_k=k$ be a solution to
Problem~\ref{pb_short1}. Let $d=(j, h)\in D$ be maximal (with
respect to the order $\succ$)
 among the elements $d_{\lambda_q}\in D$ corresponding to $q\in \{1, \ldots,
 k+1\}$. If $\Sigma'_{k+1}\subset \Om$, $\#\Sigma'_{k+1}=k+1$
  solves
Problem~\ref{pb_short1}, then $\Sigma'_{k+1}=\Sigma'_{k}\cup \{x\}$
where $x\in \Delta_p$ is in optimal position and $p\in \{1, \ldots,
k+1\}$ is such that $d=d_{\lambda_p}$ (such an element may be
nonunique).
\end{theorem}

We are able now to formulate the following statement which says
exactly how the sequence $\{\Sigma'_k\}$ of solutions to
Problem~\ref{pb_short1} looks like in the particular case we are
considering.

\begin{corollary}\label{co_div1b}
 Enumerate $D$ in the order given by the relation $\succ$, so that
 $D=\{d(i)\}_{i=1}^\infty$, $d(i+1)\succ d(i)$.
 The sequence of solutions $\Sigma'_k\subset \Om$, $\#\Sigma'_k=k$ to
Problem~\ref{pb_short1} can be described in the following way by
induction on $(D, \succ)$. Each $d(i)=(j(i), h(i))\in D$ gives rise
to a part $\mathcal{S}^i$ of the sequence of optimal sets.
\begin{itemize}
  \item $\mathcal{S}^0$ consists of the unique set $\Sigma'_1\simeq (1/2, 1/2)$.
  \item The set $\mathcal{S}^{i+1}$ consists of the sets $\Sigma'_{k+1},
  \ldots, \Sigma'_h$, where
  $\Sigma'_k$ is the last element of $\mathcal{S}^i$,
\[
\Sigma'_{j+1}=\Sigma'_{j}\cup \{x_j\}, \qquad\qquad j=k, \ldots,
h-1,
\]
$h:=k+2^{h(i)}-1$, $x_j\in \Delta_p\subset \Omega_{j(i)}$ in optimal
position, and $p\in \{1,\ldots, j\}$ is an arbitrary index
satisfying $d(i)=d_{\lambda_p}$.
\end{itemize}
\end{corollary}

\begin{proof} The proof is easily obtained by induction on the well-ordered set $(D, \succ)$.
\end{proof}

\begin{remark} An easy consequence of what has been proven so far,
is that in the case we are considering (that is, when $\nu$ is the
uniform measure over the interval) the value of $s_n$  {\em does
not\/} depend on the particular sequence of solutions
$\{\Sigma'_n\}$ to short-term Problem~\ref{pb_short1}, i.e.\ it
depends only on the index $n$.
\end{remark}

\begin{example}\label{ex_div1}
We consider the first elements of the possible sequences of
solutions $\Sigma'_k\subset \Om$, $\#\Sigma'_k=k$ to
Problem~\ref{pb_short1}.
\begin{enumerate}
  \item $i=0$, $\mathcal{S}^0=\{\Sigma'_1\}$, where
\[
\Sigma'_1\simeq (1/2, 1/2).
\]
  \item $i=1$, $d(1):=(j(1), h(1))=(-1/2, -1/2)\in D$.
  From
 $d(1)=d_{\lambda_p}$ we get $\lambda_p=3^{-j(1)+1/2}2^{-1}= 1/2$,
which corresponds to the two external intervals of length $1/2$.
Hence,  $\mathcal{S}^1:=\{\Sigma'_2,\Sigma'_3\}$, where
\begin{align*}
    \mbox{either }\Sigma'_2&\simeq (1/6, 1/3, 1/2)\mbox{ or } \Sigma'_2\simeq (1/2, 1/3,
    1/6),\\
    \Sigma'_3&\simeq (1/6, 1/3, 1/3, 1/6).
\end{align*}
  \item $i=2$, $d(2):=(j(2), h(2))=(1, 0)\in D$. From
 $d(2)=d_{\lambda_p}$ we get $\lambda_p=3^{-j(2)}2^{-h(2)}= 1/3$,
which corresponds to the two internal intervals of length $1/3$,
both belonging to $\Om_{j(2)}=\Om_1$. Hence,
$\mathcal{S}^2:=\{\Sigma'_4,\Sigma'_5\}$, where
\begin{align*}
    \mbox{either }\Sigma'_4&\simeq (1/6, 1/6,1/6, 1/3, 1/6)\mbox{ or } \Sigma'_4\simeq (1/6, 1/3, 1/6, 1/6, 1/6),\\
    \Sigma'_5&\simeq (1/6, 1/6,1/6, 1/6, 1/6, 1/6).
\end{align*}
\item $i=3$, $d(3):=(j(3), h(3))=(1/2, -1/2)\in D$.
 From
 $d(3)=d_{\lambda_p}$ we get $\lambda_p=3^{-j(3)-1/2}2^{-1}= 1/6$,
which corresponds to the two external intervals of length $1/6$.
Hence,  $\mathcal{S}^3:=\{\Sigma'_6,\Sigma'_7\}$, where
\begin{align*}
    \mbox{either }\Sigma'_4&\simeq (1/18, 1/9, 1/6,1/6, 1/6, 1/6, 1/6)\\
    \mbox{ or } \Sigma'_4& \simeq (1/6, 1/6,1/6, 1/6, 1/6, 1/9, 1/18),\\
    \Sigma'_5&\simeq (1/18, 1/9, 1/6,1/6, 1/6, 1/6, 1/9, 1/18).
\end{align*}
\item $i=4$, $d(4):=(j(4), h(4))=(1, 1)\in D$.
 From
 $d(4)=d_{\lambda_p}$ we get $\lambda_p=3^{-j(4)}2^{-h(4)}= 1/6$,
which corresponds to the four internal intervals of length $1/6$,
all belonging to $\Om_{j(4)}=\Om_1$. Hence,
$\mathcal{S}^3:=\{\Sigma'_8,\Sigma'_9,\Sigma'_{10}, \Sigma'_{11}\}$,
where
\begin{align*}
    \mbox{either }\Sigma'_8&\simeq (1/18, 1/9, 1/12, 1/12,1/6, 1/6, 1/6, 1/9, 1/18)\\
    \mbox{ or } \Sigma'_8& \simeq (1/18, 1/9, 1/6,1/12, 1/12, 1/6, 1/6, 1/9, 1/18)\\
    \mbox{ or } \Sigma'_8& \simeq (1/18, 1/9, 1/6,1/6, 1/12, 1/12, 1/6, 1/9, 1/18)\\
    \mbox{ or } \Sigma'_8& \simeq (1/18, 1/9, 1/6,1/6, 1/6, 1/12, 1/12, 1/9, 1/18),\\
    \mbox{either }\Sigma'_9&\simeq (1/18, 1/9, 1/12, 1/12,1/12, 1/12, 1/6, 1/6, 1/9, 1/18)\\
        \mbox{or }\Sigma'_9&\simeq (1/18, 1/9, 1/12, 1/12,1/6, 1/12, 1/12,  1/6, 1/9, 1/18)\\
   \mbox{or }\Sigma'_9&\simeq (1/18, 1/9, 1/12, 1/12,1/6,1/6, 1/12, 1/12,  1/9, 1/18)\\
 \mbox{ or } \Sigma'_9& \simeq (1/18, 1/9, 1/6, 1/12, 1/12,  1/12, 1/12,1/6, 1/9, 1/18)\\
 \mbox{ or } \Sigma'_9& \simeq (1/18, 1/9, 1/6, 1/12, 1/12,1/6,  1/12, 1/12, 1/9, 1/18)\\
 \mbox{ or } \Sigma'_9& \simeq (1/18, 1/9, 1/6, 1/6, 1/12, 1/12,  1/12, 1/12, 1/9, 1/18),\\
 & \ldots\ldots
\end{align*}
\end{enumerate}
\end{example}

We are able to claim now the following assertions answering
questions~(B) and~(C) for the case we are considering.

\begin{theorem}\label{th_limnonunique}
Assume $\nu=\lcal^1\res [0,1]$. Then the bounded sequence $\{n
s_n\}$ does not converge to a limit as $n\to \infty$. Further, for
any sequence $\{\Sigma'_n\}$ of solutions to
Problem~\ref{pb_short1}, the sequence of probability measures
$\mu_{\Sigma'_n}$ has infinitely many limit measures in the weak$^*$
sense as $n\to \infty$.
\end{theorem}

\begin{proof}
Fixed a $j\in \N$,  consider the sequence of indices $k_n$ such that
$\Sigma'_{k_n}$ induces a subdivision of the set $\Omega_j$ into
$2^n$ equal subintervals while the partition of $\Om$ corresponding
to $\Sigma'_{k_n+1}$ divides further one of the latter subintervals
into two equal intervals. This means that the sequence of partitions
corresponding to optimal sets $\Sigma'_k$, starting from $k=k_n$ and
up to $k=k_n+2^n$, will be obtained by dividing into two equal parts
at each step one of the $2^n$ subintervals in $\Omega_j$. Set now
$$
a_n:=\int_{[0,1]\setminus
\Omega_j}\!\!\dist(x,\Sigma'_{k_n})dx\quad\mbox{ and }\quad
b_n:=\int_{\Omega_j}\dist(x,\Sigma'_{k_n})dx.
$$
Fix an arbitrary dyadic number $\lambda\in (0,1)$ and consider the
new sequence of indices $k'_n:=k_n+\lambda 2^n$, so that the
partition of $\Om$ corresponding to $\Sigma_{k'_n}$ coincides with
that corresponding to $\Sigma_{k_n}$ up to the fact that some of the
subintervals in $\Omega_0$, namely, a fraction $\lambda$ of the
total, have been split in two parts. This implies that
\[
k'_nl_{k'_n}=k'_nF(\Sigma'_{k'_n})=(k_n+\lambda
2^n)(a_n+(1-\lambda/2)b_n),
\]
splitting an interval in two equal parts reduces the average
distance from its points to $\Sigma$  by a factor two. Mind that
\[
k_nl_{k_n}=k_nF(\Sigma'_{k_n})=k_n(a_n+b_n),
\]
and assume  by contradiction that, $\lim_k k l_k$ exists. Then
$\lim_n k'_nl_{k'_n} = \lim_n k_nl_{k_n}$, and hence,
\begin{align*}
(k_n+\lambda 2^n)\Big(a_n+ &
\Big(1-\frac{\lambda}{2}\Big)b_n\Big)-k_n(a_n+b_n)\\
&=\lambda\Big(2^na_n-\frac{1}{2}k_nb_n+\left(1-\frac{\lambda}{2}\right)2^nb_n\Big)\\
&=\lambda\Big(2^na_n-\frac{1}{2}k_nb_n+\frac{1}{9}\Big(1-\frac{\lambda}{2}\Big)\Big)\to
0
\end{align*}
as $n\to \infty$.
 Therefore, $2^na_n-k_nb_n/2\to \lambda/18-1/9$ as $n\to \infty$,
 which means that, taking two different dyadic values of $\lambda\in (0,1)$,
 the sequence $\{2^na_n-k_nb_n/2\}$ has two different limits. This
 contradiction proves the first claim.

To prove the second claim, we first show that the limit measures of
the sequence $\{\mu_{\Sigma'_{l}}\}$ are not unique, where $l\in
[k_n, k_n+2^n]$. To this aim, suppose by contradiction that all the
subsequences of the above sequence converge in the weak$^*$ sense to
the same limit measure $\mu$ as $n\to \infty$. In particular, this
is the case of $\{\mu_{\Sigma'_{k_n}}\}$, which clearly then
converges to a uniform measure over $\Omega_j$ since the points of
$\Sigma'_{k_n}$ are uniformly distributed over $\Omega_j$. Observing
that $k_n\geq 2^n$, and hence
\[
1/(1+\lambda) \leq k_n/k_n' \leq 1,
\]
we may assume without loss of generality that up to extracting a
subsequence of $k_n$ (not relabeled), the sequence $\{k_n/k_n'\}$
converges to some finite limit. Hence
\[
1= \lim_{n\to
\infty}\frac{\mu_{\Sigma'_{k_n'}}(\Omega_j)}{\mu_{\Sigma'_{k_n}}(\Omega_j)}=
\lim_{n\to \infty}\frac{(1+\lambda) 2^n/ k_n'}{2^n/k_n}= (1+\lambda)
\lim_{n\to \infty}\frac{k_n}{k_n'},
\]
which gives $\lim_n k_n/k_n' \neq 1$. But since
\[
\mu_{\Sigma'_{k_n'}}(\Omega_j)
=\frac{k_n}{k_n'}\mu_{\Sigma'_{k_n}}(\Omega_j)
 + \frac{\lambda 2^n}{k_n'},
\]
then passing to a limit in the above relationship as $n\to \infty$,
we get
\begin{align*}
\mu(\Omega_j) &=\lim_{n\to \infty}\frac{k_n}{k_n'}\mu(\Omega_j)
 + \lim_{n\to \infty}\frac{\lambda 2^n}{k_n'}\\
 & = \mu(\Omega_j) \lim_{n\to \infty}\frac{k_n}{k_n'}
 + \lim_{n\to \infty} \left(1-\frac{k_n}{k_n'}\right).
\end{align*} Minding that $\lim_n k_n/k_n' \neq 1$, the above
relationship is only possible when $\mu(\Omega_j)=1$, i.e.\ when the
unique limit measure $\mu$ is concentrated on $\Omega_j$, which is
clearly not
the case because
\begin{equation}\label{eq_mjplus1}
    \mu(\Omega_{j+1}) \geq \frac {\mu(\Omega_j)}{ 3}.
\end{equation}
To verify the latter inequality, it is enough to notice that
\[
\frac{\mu_{\Sigma'_{k_n}}(\Omega_{j+1})}{\mu_{\Sigma'_{k_n}}(\Omega_{j})}=
\frac{2^{h}+1}{ 2^n+1}, \mbox{ with } (j,n) \succ (j+1, h),
\]
so that $h\geq n-\log 3/\log 2$, which gives~\eqref{eq_mjplus1} in
the limit as $n\to \infty$.


At last, to show that the limit measures of the sequence
$\{\mu_{\Sigma'_n}\}$ are infinite, it is enough to vary $j\in \N$.
\end{proof}

\begin{figure}[h]
\includegraphics[height=6cm]{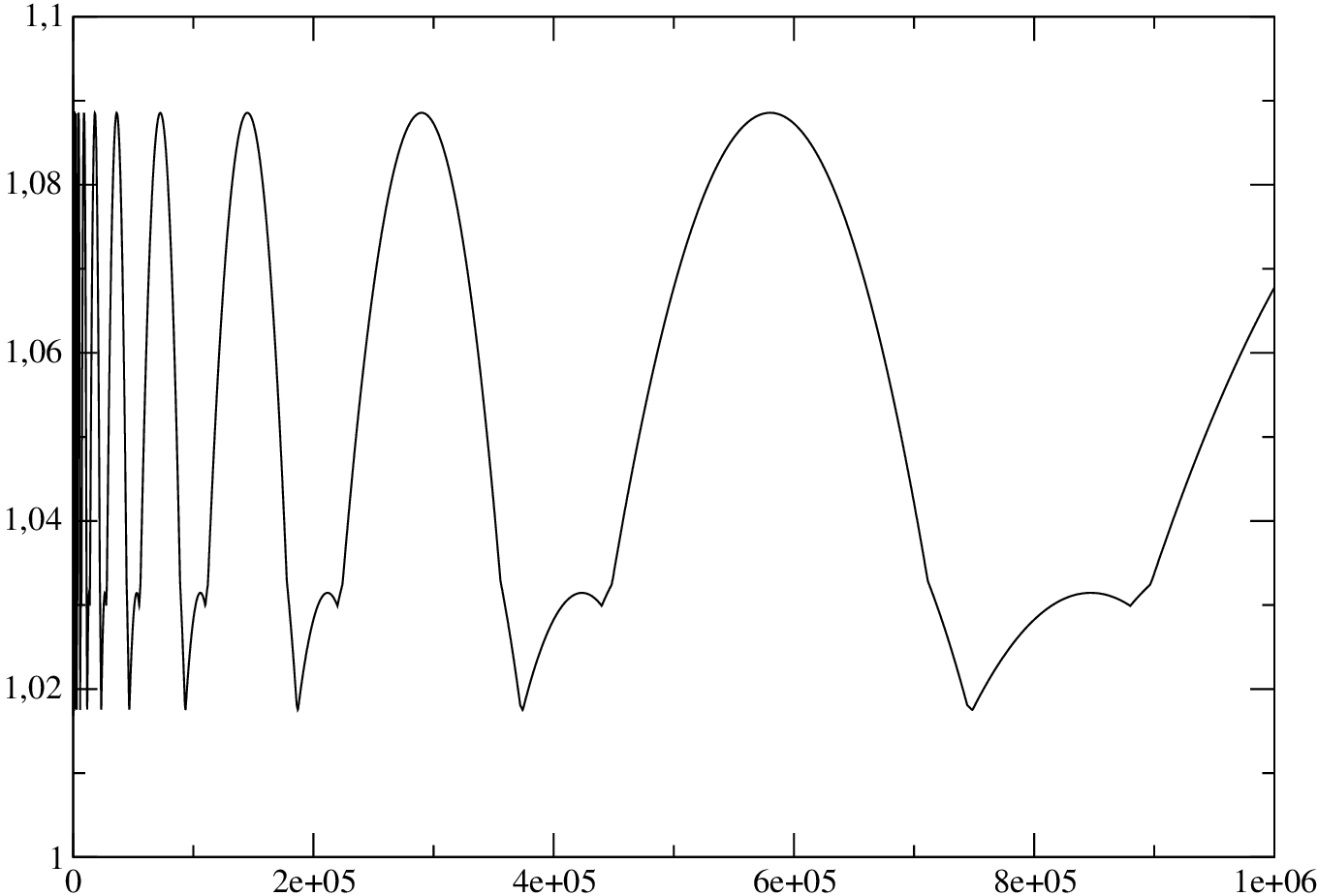}
\caption{Plot of $s_n/l_n$ for the case of uniform density over
$[0,1]$.}\label{SLP/LLP}
\end{figure}

In Figure~\ref{SLP/LLP} below we provide a graph of $s_n/l_n$ for
the case we are considering calculated according to the algorithm
given by Theorem~\ref{th_div1} and~Corollary~\ref{co_div1b}.

\vspace{1cm}

\begin{theorem}\label{thm:lebesgue_is_not_a_limit}
No sequence $\{\Sigma'_n\}$ of solutions to Problem~\ref{pb_short1}
(with $\nu$ uniform measure on $[0,1]$) has the Lebesgue measure
$\lcal^1 \res [0,1]$ as a limit measure of some subsequence of
$\mu_{\Sigma'_n}$ in the weak$^*$ sense as $n\to \infty$.
\end{theorem}

\begin{proof}
By using Corollary~\ref{co_div1b}, to any $k$ we can associate a set
$\Omega_j$ in the following way. Let $i\in \N$ be such that
$\Sigma'_k\in \mathcal{S}^i$. Consider now the pair $(j(i),h(i))\in
D$ corresponding to the index $i$, and set $j:=j(i)$. The index
$j=j[k]$ associated to a set $\Sigma'_k$ represents the set
$\Omega_j$ where we last put a point in building $\Sigma'_k$ and
also the set $\Omega_j$ which is being divided in a dyadic way by
the optimal sequence at step $k$. Namely, if we look at the points
composing $\Sigma'_k$, we see the following picture.
\begin{itemize}
\item For every $j> \max_{k'\leq k} j[k']$ no point of $\Sigma_k$ belongs to
$\Omega_j$.
\item Every $\Omega_j$ with $j\leq \max_{k'\leq k} j[k']$ and $j\neq j[k]$ contains a certain number
$2^h$ of points of $\Sigma'_k$ (precisely we have $h=\max\{h\,:\,
(j,h)\prec(j[k],k)\}$).
\item At last, in $\Omega_{j[k]}$ we cannot exactly predict the number of points of $\Sigma'_k$.
\end{itemize}
Up to choosing a subsequence of $n$ (not relabeled), we can suppose
the existence of two distinct indices $j_1$ and $j_2$ such that for
any $n$ (from the chosen subsequence) we have $j[n]\neq j_1$, and
$j[n]\neq j_2$. Extracting a further subsequence of $n$ (again not
relabeled), we may assume that
 $\max_{k'\leq n} j[k']\geq (j_1\vee j_2)$, and hence
the number of points of $\Sigma'_{n}$ in the sets $\Omega_{j_1}$ and
$\Omega_{j_2}$ is a power of two. This means that, if we set
$\mu_n:=\mu_{\Sigma'_{k_n}}$ and we suppose
$\mu_n\destar\mu=\lcal^1\res [0,1]$ in the weak$^*$ sense as
$n\to\infty$, we get
$$
2^{\lambda(n,j_1,j_2)}=
\frac{\mu_n(\Omega_{j_1})}{\mu_n(\Omega_{j_2})}\to
\frac{\mu(\Omega_{j_1})}{\mu(\Omega_{j_2})}=3^{j_2-j_1},
$$
where $\lambda(n,j_1,j_2)$ is an integer exponent depending only on
$n$, $j_1$ and $j_2$. Yet this is a contradiction as a sequence of
powers of two can converge only to $0$, $+\infty$ or a power of two.
\end{proof}

The above theorem gives a rigorous formulation of the insight of the
authors of~\cite{TSuzAsamOkab91} who noticed that a short-term
strategy in the above one-dimensional situation, compared to other
allocation policies, induces ``a more uniformly spaced allocation
yet not {\em completely uniform}''.

\begin{corollary}
One has
\[
  1<\liminf_{n\to \infty} s_n/l_n <\limsup_{n\to \infty} s_n/l_n,
  \]
  if $\nu$ is the uniform measure over $[0,1]$.
  \end{corollary}

\begin{proof}
The second inequality is just the reformulation of the first claim
of Theorem~\ref{th_limnonunique} minding
Corollary~\ref{consequences}. To prove the first inequality, note
that $s_n\geq l_n$, while,
 should the above $\liminf$ be
equal to one, we could build a subsequence $\{\Sigma'_{k_n}\}$ of
solutions to short-term Problem~\ref{pb_short1} such that the
respective sequence of measures $\{\mu_{\Sigma'_{k_n}}\}$ is
asymptotically optimal for the sequence of functionals
$\{\F_{k_n}\}_n$ defined by~\eqref{eq_FnProb}. From general results
in $\Gamma-$convergence theory we would obtain
$\mu_{\Sigma_{k_n}}\destar\bar{\mu}$ in the weak$^*$ sense as
$n\to\infty$, where $\bar{\mu}$ is the limit measure for the
long-term Problem~\ref{pb_long1}, i.e. the Lebesgue measure over
$[0,1]$, which  contradicts
Theorem~\ref{thm:lebesgue_is_not_a_limit}.
\end{proof}

\section{Concluding remarks and open problems}\label{lastsec}

We conclude the paper by a list of remarks and open problems
regarding the model studied in this paper, which we consider to be
interesting for further study.

\bigskip
{\bf Question I}. 
  It seems
  interesting to understand whether the phenomena
studied in Section~\ref{onedim} (e.g.\ non existence of the limit of
$s_n/l_n$, the fact that the respective $\liminf$ is strictly
greater than one, non-uniqueness of the limit measures for the
short-term location problem etc.) for the case $\nu=\lcal^1\res
[0,1]$
  occur in the case of generic measure
  $\nu$ in one-dimensional case $d=1$. Further, it seems to be important for applications to
  obtain sharp estimates on $\liminf_{n\to \infty} s_n/l_n$ and $\limsup_{n\to \infty}
  s_n/l_n$ even for the uniform measure $\nu$ (in the latter case,
  the numerical computation provided in this paper (Figure \ref{SLP/LLP}) suggests that
  both limits are close to $1$, so that, in practice, one could assume that
  asymptotically the value of the functional in the short-term problem is almost the same
  as in the long-term one, which would give a rigorous statement of the idea
first suggested in~\cite{TSuzAsamOkab91}) and to characterize
completely the limit measures of the short-term problem.

\begin{figure}[h]
\includegraphics[width=8cm]{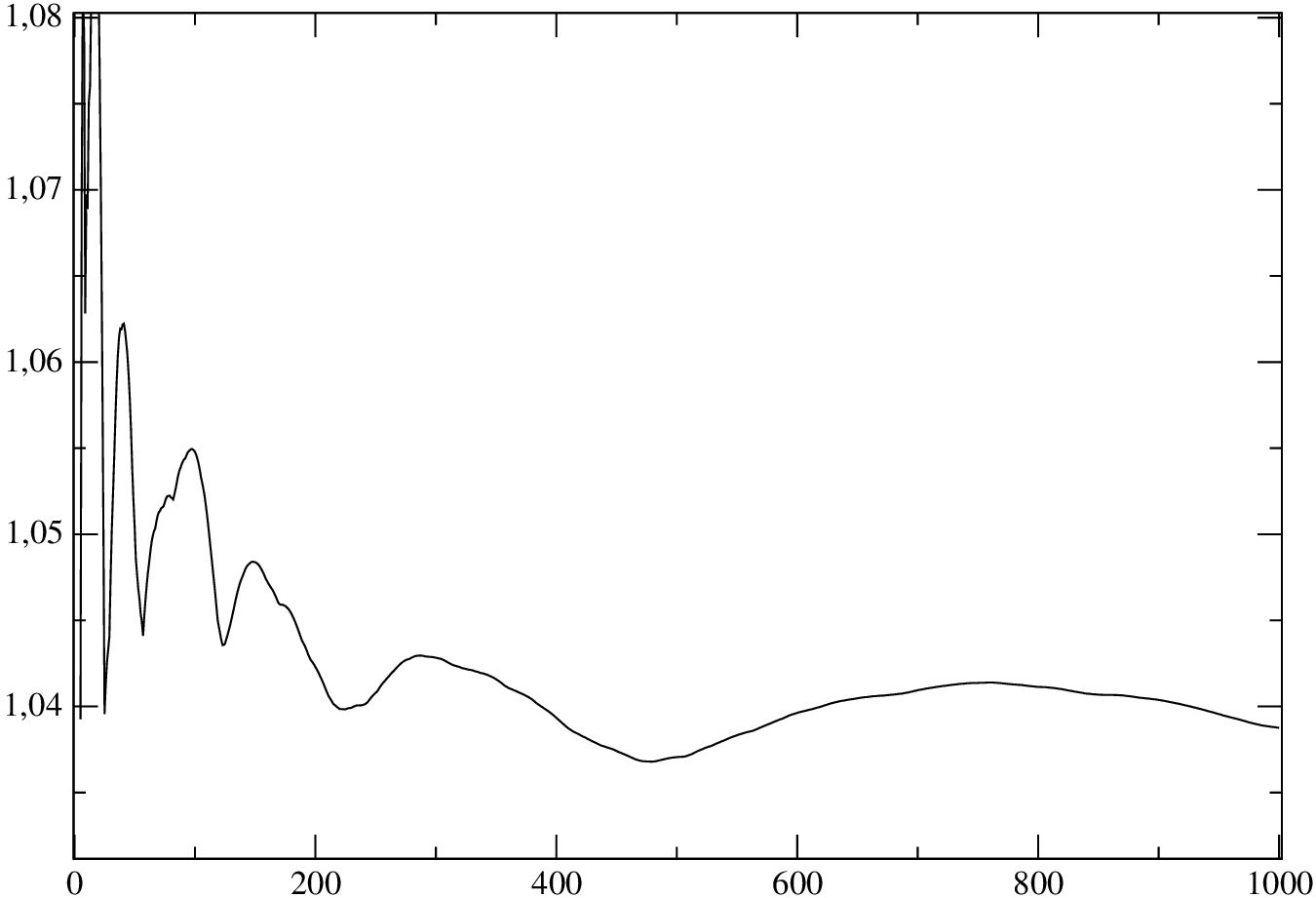}
\caption{Plot of $s_n/l_n^\infty$  
for the case of a uniform density on a unit
square}\label{lp1000p200g-graph}
\end{figure}
\begin{figure}[h]
\includegraphics[width=6cm]{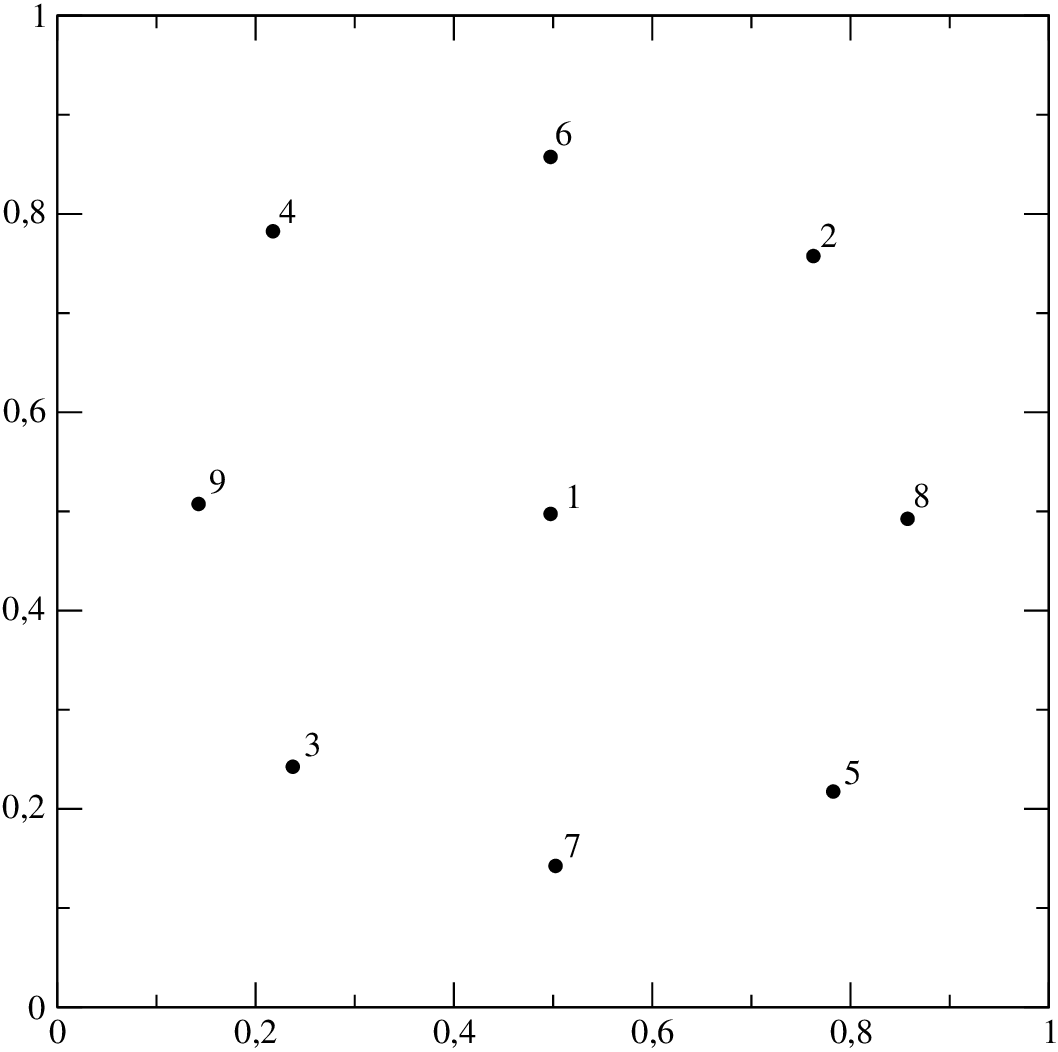}
\hspace{1ex}
\includegraphics[width=6cm]{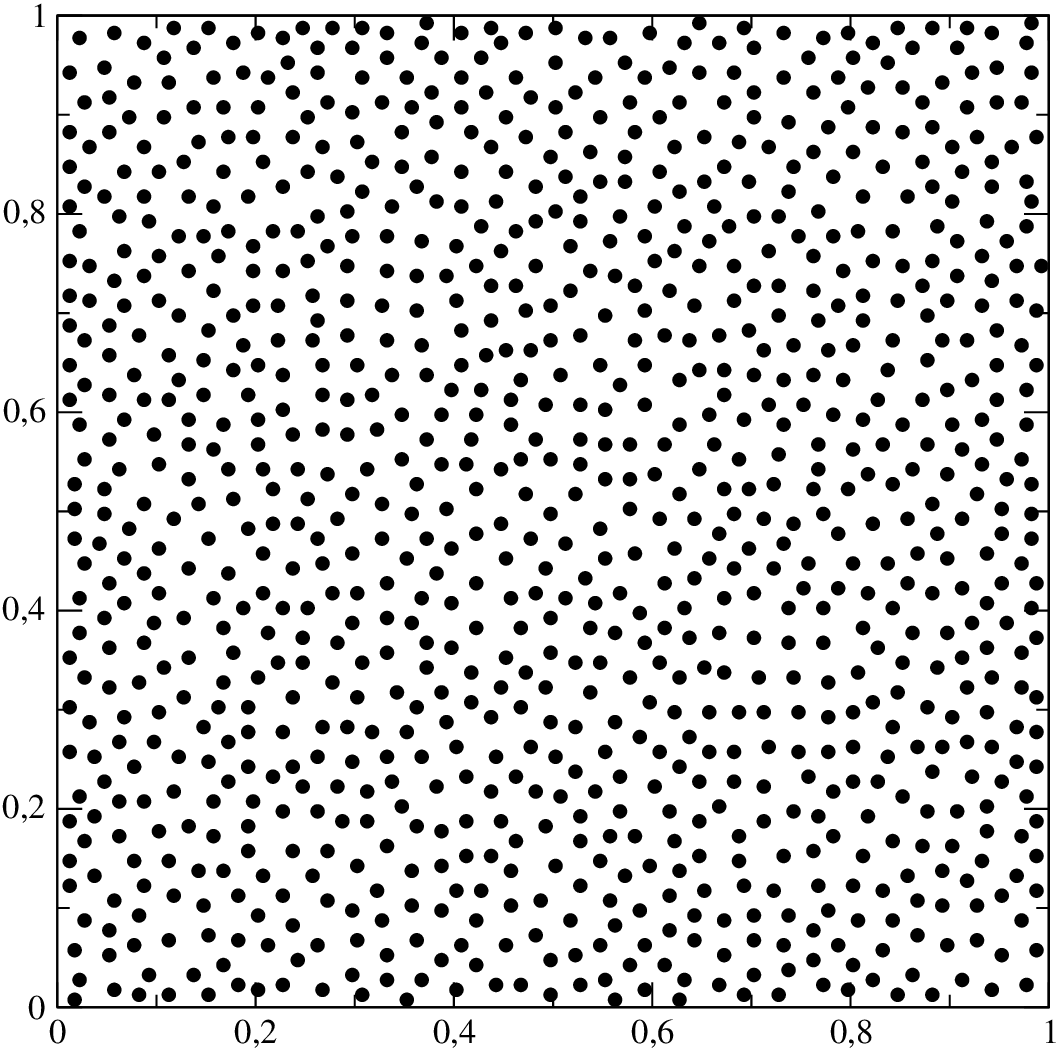}
\caption{Short-term location for the case of a uniform density on a
unit square of (left) the first $9$ points (right) $1000$ points.
}\label{lp1000p200g-set}
\end{figure}

\bigskip
  {\bf Question II}. In the case of the generic space dimension $d >1$ we
  are able to prove only the asymptotic order estimate of the value of the functional in the short-term
  problem. It seems therefore important to study the above posed
  problems for generic space dimension even in the simplest case
  when, say, $\nu$ is the uniform measure on the unit ball or on the unit square.
  Figure~\ref{lp1000p200g-graph} shows what one can expect
  about the behavior of $s_n$ in the case of a uniform density on a unit
square (for convenience, the normalized values $s_n/l_n^\infty$ are
provided, where $l_n^\infty:=n^{-1/d}\theta_2 \|f\|_{d/(d+1)}$
stands for the asymptotical value of  the long-term minima
according to Theorem~\ref{th_UCGamloc}), while
Figure~\ref{lp1000p200g-set} shows how the distribution of
  points in this case looks like. Both figures are obtained by a
numerical calculation on a uniform $200\times 200$ grid. More
results of numerical computations for the short-term optimal
location problem can be fund on the web
page~\cite{Web-page-brabutsan06}.

\vspace{1cm}


\bigskip
{\bf Question III}. At last, it is natural to mention here a similar
  problem introduced in~\cite{ButOudSte02} and sometimes called {\em
  irrigation problem\/} (see~\cite{BrancoliniButtazzo,ButSte03,ButSte04pal,PaoSte04,TilliSantAmb04,Ste05a})
  on minimization of the average distance functional but over  compact connected sets of finite length rather than
  over discrete sets of points like in this paper.
  The statement of such a problem is obtained by
  replacing the constraint on cardinality $\#\Sigma$ of the unknown minimizer $\Sigma$
  in the location
problem  by the similar constraint on the one-dimensional Hausdorff
measure $\HH^1(\Sigma)$. The problem is therefore that of finding a
minimizer of the cost $F\colon \Sigma\mapsto
\int_{\Omega}\dist(x,\Sigma)\,d\nu(x)$ over all compact and
connected sets $\Sigma\subset\bar{\Omega}$ satisfying
$\HH^1(\Sigma)\leq l$ with given $l>0$. This problem, which again
can be interpreted as the long-term one, admits also a natural
short-term approach. Namely, by letting the parameter $l$ increase,
we would like to find an irrigation set $\Sigma$ which increases
continuously with $l$, in such a way that the functional $F$
decreases as fast as possible. This construction may be made
rigorous through a slight modification of the well-known method of
{\em minimizing movements\/} (see for
instance~\cite{Ambrosio,AmbrosioGigliSavare} for the presentation of
the theory). Namely, fixed an arbitrary time step $\tau>0$, we
minimize the functional $\Sigma\mapsto
F(\Sigma)+\HH^1(\Sigma\setminus\Sigma_k^\tau)^2/2\tau$ over the set
of all connected compact subsets of $\overline{\Omega}$ satisfying
the constraint $\Sigma_{k+1}^\tau\supset\Sigma_k^\tau$, where
$\Sigma_0^\tau:=\emptyset$ (alternatively, one could drop the
$\HH^1-$penalization term in the above functional and add the
additional constraint $\HH^1(\Sigma_{k+1})\leq
\huno(\Sigma_{k}^\tau)+\tau$ instead). Let then $t\mapsto
\Sigma^\tau(t)$ be the piecewise constant map defined by
$$
\Sigma^\tau(t):=\Sigma^\tau_k\qquad\hbox{for }t\in[k\tau,(k+1)\tau[.
$$
One should first study whether in this way one can find a
well-defined increasing evolution $\Sigma(t)$ as a limit (in the
Hausdorff topology)  of $\Sigma^\tau(t)$ as $\tau\to 0$. Further,
several questions then arise about the behavior of $\Sigma(t)$. For
instance, it has been proven
in~\cite{ButSte03,ButSte04pal,TilliSantAmb04,Ste05a} that the
solutions of the long-term irrigation problem under suitable
conditions on problem data do not contain loops, have a finite
number of endpoints and may have only regular tripods (i.e.\ triple
junctions with the branches with infinitesimal angles $120^\circ$
between each other) as branching points, which are at most finite in
number. It is interesting to verify whether the same or similar
properties hold also for solutions of the short-term problem
$\Sigma(t)$ for all (or for some) $t\in \R^+$, or whether such
solutions are just simple curves without branching points. Besides,
one is also interested in the asymptotic density of the curve
$\Sigma(t)$ (i.e.\ the weak$^*$ limits of the measures
$\HH^1\res\Sigma(t)/\HH^1(\Sigma(t))$ as $t\to \infty$), which for
the long-term irrigation problem was studied in~\cite{TillMosc02}.

%
%

\bigskip
{\small
\begin{minipage}[t]{5.7cm}
Alessio Brancolini\\
SISSA\\
Via Beirut, 4\\
34014 Trieste, ITALY\\
{\tt brancoli@sissa.it}
\end{minipage}
\begin{minipage}[t]{5.7cm}
Filippo Santambrogio\\
Scuola Normale Superiore\\
Piazza dei Cavalieri, 7\\
56126 Pisa, ITALY\\
{\tt f.santambrogio@sns.it}
\end{minipage}

\bigskip

\begin{minipage}[t]{5.7cm}
Giuseppe Buttazzo\\
Dipartimento di Matematica\\
Universit\`a di Pisa\\
Largo B. Pontecorvo, 5\\
56127 Pisa, ITALY\\
{\tt buttazzo@dm.unipi.it}
\end{minipage}
\begin{minipage}[t]{5.7cm}
Eugene Stepanov\\
Dipartimento di Matematica\\
Universit\`a di Pisa\\
Largo B. Pontecorvo, 5\\
56127 Pisa, ITALY\\
{\tt e.stepanov@sns.it}
\end{minipage}
}

\end{document}